\documentclass[11pt]{article}
\usepackage{amsfonts}
\usepackage{latexsym,amsmath}
\usepackage{amssymb,array}
\makeatletter \oddsidemargin 0in \evensidemargin 0in \textwidth
16cm \RequirePackage[dvips]{graphicx} \textheight 20cm
\newtheorem{t1}{Theorem}[section]

\newtheorem{d1}{Definition}[section]
\newtheorem{r1}{Remark}[section]

\begin{document}
\title{Inequalities involving expectations of selected functions in reliability theory to characterize distributions}
\author{{\bf Chanchal Kundu}\footnote{Corresponding
author e-mail: ckundu@rgipt.ac.in;
chanchal$_{-}$kundu@yahoo.com.}\\\and {\bf Amit Ghosh}
\and\and Department of Mathematics\\
Rajiv Gandhi Institute of Petroleum Technology\\
Rae Bareli 229 316, U.P., India}
\date{Revised version to appear in {\it Communications in Statistics$-$Theory \& Methods},\\ $\copyright$ by Taylor $\&$ Francis Group.\\Submitted: July, 2015}
\maketitle
\begin{abstract}
Recently, authors have studied inequalities involving expectations of selected functions viz. failure rate, mean residual life, aging intensity function and log-odds rate which are defined for left truncated random variables in reliability theory to characterize some well-known distributions. However, there has been growing interest in the study of these functions in reversed time and their applications. In the present work we consider reversed hazard rate, expected inactivity time and reversed aging intensity function to deal with right truncated random variables and characterize a few statistical distributions.
\end{abstract}
\section{Introduction}
In the literature, the problem of characterizing probability distributions has been investigated by many researchers. In order to fit a particular probability distribution to the real world data, it is essential to confirm whether the given probability distribution satisfying the underlying requirements by its characterization. In fact, characterization approach is very appealing to both theoreticians and applied workers. A probability distribution can be characterized through various methods. In recent years, authors have considered the approach of studying inequalities involving expectations of various functions in reliability theory to characterize distributions. For example, Nanda (2010) studied the characterizations of distributions through the expected values of failure rate (FR) and mean residual life (MRL) functions. Bhattacharjee et al. (2013) consider FR, MRL, log-odds rate and aging intensity (AI) functions to characterize a few distributions.\\
\hspace*{.2in} All the functions they have considered deal with left truncated random variables. But several common situation may occur when data are left-censored or right-truncated. For example, in the study of inequality of income distribution it is of interest to study the inequality of a population eliminating high (richest population) values.  Analogous measures say, reversed hazard rate (RHR), expected inactivity time (EIT) etc. are defined in the literature to deal with right truncated distribution. These functions characterize the aging phenomenon of any living unit or a system of components which has been found to be failed at time $t$. Let $X$ be an absolutely continuous random variable with probability density function  $f$ and cumulative distribution function $F$. Let the support of the
random variable $X$ be $(a,b)$ where $a=\inf\{t:F(t)>0\}$ and $b=\sup\{t:F(t)<1\}$ with $-\infty\leqslant a<b<\infty$. If $a$ is taken to be zero, the random variable $X$ may be thought of as the random lifetime of a living organism or of a system of components. Let $X$ have a finite second-order moment given by $Var(X)=\sigma^2$ with mean $E(X)=\mu$. For any positive integer $r$ let us denote the $r^{th}$ order raw moment of $X$ as $\mu_r$ with $\mu_{1}=\mu$, where there is no ambiguity.\\
\hspace*{.2in} The RHR function of $X$ is defined as $\phi(t)=f(t)/ F(t)$, where defined. If a system is down at time $t$, then the probability that it was active at time $t-\epsilon$, for $\epsilon (>0)$ very small is given approximately by $\epsilon\phi(t)$. With this physical significance, the inactivity time is defined as the random variable $X_t=(t-X|X\leqslant t)$ where the EIT, also called mean inactivity time (MIT), is given by \begin{eqnarray}\label{eq0}m(t)=E(X_t)=\frac{1}{F(t)}\int_a^tF(x)dx.
\end{eqnarray}
The RHR and EIT are widely used in many areas of applied probability and statistics such as forensic study, actuarial science, biometry, survival analysis, economics, risk management and reliability. They are closely related in the sense that both of them determine a distribution uniquely via the relation
\begin{eqnarray}\label{eq1}F(t)=\exp\left[-\int_t^b \phi(u)du\right],
\end{eqnarray}
where $\phi(t) m(t) =1-m'(t)$. For characterization of distribution through RHR and EIT one may refer to Chandra and Roy (2001), Kundu and Nanda (2010), Kundu et al. (2010) and Kundu and Sarkar (2015), among others. Different properties of EIT have been discussed in Kayid and Ahmad (2004), Ahmad et al. (2005), Li and Xu (2006), Zhang and Cheng (2010) and Kayid and Izadkhah (2014). In analogy with Jiang et al. (2003), the reversed AI function of $X$ is defined as follows.
\begin{d1} For an absolutely continuous random variable $X$ having support $(a,b)$, the reversed AI function of $X$ at $t$ is defined as
\begin{eqnarray*}
\overline L(t)&=&\frac{(b-t)\phi(t)}{\int_t^b\phi(u)du},~a<t<b\\
&=&\frac{(t-b)f(t)}{F(t)\ln F(t)}.
\end{eqnarray*}
\end{d1}
This may be considered as a dual measure of AI function to evaluate the aging property of a unit (down at time $t$) quantitatively.
\begin{r1} It is to be noted that $\overline L(t)$ is not defined for $t\leqslant a$ whereas this may be defined as zero when $t\geqslant b$.
\end{r1}
\hspace*{.2in} It can be easily verified that the reversed AI function $\overline L(t)=1$ if and only if the RHR function $\phi(t)=\gamma$ (a positive constant), which characterizes Type 3 extreme value distribution having cumulative distribution function
\begin{equation}\label{eq2}F(t)=\left\{\begin{array}{ll}e^{\gamma(t-b)},&
t\in(-\infty,b],~b\geqslant0\\
 1,& otherwise.\end{array}\right.
 \end{equation}
Further, $\overline L(t)<1$ if RHR is increasing and $\overline L(t)>1$ if RHR is decreasing. The smaller the value of $\overline L(t)$, the stronger the tendency of aging of the random variable $X$.\\
\hspace*{.2in} Looking into the importance, here we consider RHR, EIT and reversed AI function to characterize a few distributions. The rest of the paper is arranged as follows. In Section 2 we provide characterizations of quite a few useful continuous
distributions through the expected values of RHR. Inequalities involving expectation of EIT is considered in Section
3 along with some related study based on reversed AI function.
\section{Inequalities involving expectation of RHR}
This section characterizes some well-known distributions through the expected values of RHR. These characterization results may be helpful for model selection in reliability theory and survival analysis.\\
\hspace*{.2in} In the following theorem we characterize Type 3 extreme value distribution which also belongs to the reversed generalized Pareto distribution developed by Castillo and Hadi (1995) as a fatigue model that satisfy certain compatibility conditions arising out of physical and statistical conditions in fatigue studies.
\begin{t1}\label{th2.1}
For any absolutely continuous random variable $X$,
$$E\left[\frac{1}{\phi(X)}\right]\geqslant\frac{1}{E[\phi(X)]}.$$
The equality holds if and only if $X$ follows the distribution given in (\ref{eq2}) for a suitable choice of $\gamma$.
\end{t1}
Proof: By Cauchy-Schwartz inequality, we have
\begin{equation} \label{eq2.2}E\left[\frac{1}{\phi(X)}\right]E[\phi(X)]\geqslant~1 .
\end{equation}
The equality in (\ref{eq2.2}) holds if and only if there exists a constant $A(>0)$ such that, for all $x\in(-\infty,b]$,
$$\frac{F(x)}{f(x)}f(x)=A\frac{f(x)}{F(x)}f(x),$$
which is equivalent to the fact that $\phi(x)=1/\sqrt{A}=\gamma$, say. Now, on using (\ref{eq1}), we get the required result. $\hfill\square$\\

\hspace*{.2in} Below we characterize power distribution.
\begin{t1}\label{th2.2}
Let $X$ be an absolutely continuous nonnegative random variable with $E[X\phi(X)]<\infty$ and $E\left[\frac{1}{X\phi(X)}\right]<\infty $, then
\begin{eqnarray}\label{eq2.3}E\left[\frac{1} {X \phi(X)}\right]\geqslant\frac{1}{E\left[X\phi(X)\right]},
\end{eqnarray}
and the equality holds if and only if $X$ follows power distribution
\begin{equation*} F(x)=\left\{\begin{array}{ll}\left(\frac{x}{b}\right)^{c},&
0<x<b,\;b,c>0\\
 1,& otherwise.\end{array}\right.
\end{equation*}
\end{t1}
Proof: By Cauchy-Schwartz inequality, we obtain (\ref{eq2.3}). The equality in (\ref{eq2.3}) holds if and only if there exists a constant $A'(>0)$ such that
$$\frac{F(x)}{xf(x)}f(x)=A'\frac{xf(x)}{F(x)}f(x),$$
which is equivalent to the fact that $\phi(x)=c/x,$ where $c=1/\sqrt{A'}$ is a constant. Now, on using (\ref{eq1}), we get the required result. $\hfill\square$\\

\hspace*{.2in} In the following theorem a characterization of inverse Weibull distribution (cf. Keller and Kamath, 1982) also known as complementary Weibull or reciprocal Weibull distribution is provided.
\begin{t1}\label{th2.4}
Let $X$ be an absolutely continuous nonnegative random variable with $E[X^k\phi(X)]<\infty$ and $E\left[\frac{1}{X^k\phi(X)}\right]<\infty$. Then, for some positive integer $k(>1),$
\begin{eqnarray}\label{eq2.5}E\left[\frac{1} {X^k \phi(X)}\right]\geqslant \frac {1} {E\left[X^k\phi(X)\right]},
\end{eqnarray}
and the equality holds if and only if $X$ follows inverse Weibull distribution given by
$$ F(x)=e^{-\nu x^{-\delta}}; \nu>0,\delta>0,~x\in(0,\infty).$$
\end{t1}
Proof: By Cauchy-Schwartz inequality, we obtain (\ref{eq2.5}).
The equality in (\ref{eq2.5}) holds if and only if there exists a constant $A^*(>0)$ such that, for all $x\in(0,\infty)$,
$$\frac{f(x)}{x^k\phi(x)}=A^*{x^k\phi(x)f(x)},$$
which is equivalent to the fact that $\phi(x)=\theta x^{-k},$ where $\theta>0$ is a constant. Now, on using (\ref{eq1}), we get
$$ F(x)=e^{\frac{\theta}{k-1}\left(-\frac{1}{x^{k-1}}\right)} ~,~x\in(0,\infty).$$
Hence the result follows  by noting that $\nu(=\theta/(k-1))>0$ and $\delta(=k-1)>0$.$\hfill\square$
\begin{r1}
The results for $k=0$ and $k=1$ have already been discussed in Theorem \ref{th2.1} and Theorem \ref{th2.2}, respectively.
\end{r1}
\hspace*{.2in} Below we characterize power function distribution.
\begin{t1}\label{th2.5}
Let $X$ be an absolutely continuous random variable with $E[e^X\phi(X)]<\infty$ and $E\left[\frac{1}{e^X\phi(X)}\right]<\infty$, then
\begin{eqnarray}\label{eq2.6}E\left[\frac{1} {e^X \phi(X)}\right]\geqslant\frac{1}{E\left[e^X\phi(X)\right]}
\end{eqnarray}
and the equality holds if and only if $X$ follows $F^\alpha$ distribution where the base line distribution is truncated extreme value distribution having cumulative distribution function
$$ F(x)=e^{-(e^{-x}-e^{-b})},~x\in(-\infty,b].$$
\end{t1}
Proof: The inequality in (\ref{eq2.6}) is a simple consequence of Cauchy-Schwartz inequality and
the equality holds if and only if there exists a constant $B(>0)$ such that, for all $x\in(-\infty,b]$,
$$\frac{f(x)}{e^x\phi(x)}=B{e^x\phi(x)f(x)},$$
which is equivalent to the fact that $\phi(x)=\alpha e^{-x}, where ~\alpha >0.$ Now, on using (\ref{eq1}), we get the required result. $\hfill\square$\\
\begin{t1}\label{th2.6}
Let $X$ be an absolutely continuous random variable with $E[a^X\phi(X)]<\infty$ and $E\left[\frac{1}{a^X\phi(X)}\right]<\infty$, where $a>1$. Then
\begin{eqnarray}\label{eq2.7}E\left[\frac{1} {a^X \phi(X)}\right]\geqslant\frac{1}{E\left[a^X\phi(X)\right]},
\end{eqnarray}
 and the equality holds if and only if $X$ follows the distribution having cumulative distribution function
 $$ F(x)=e^{-\theta(a^{-x}-a^{-b})}~;~x~\in(-\infty,b],~\theta>0.$$
\end{t1}
Proof: By Cauchy-Schwartz inequality, we obtain (\ref{eq2.7}).
The equality in (\ref{eq2.7}) holds if and only if there exists a constant $B'(>0)$ such that, for all $x\in(-\infty,b]$,
$$\frac{f(x)}{a^x\phi(x)}=B'{a^x\phi(x)f(x)},$$
which is equivalent to the fact that $\phi(x)=\theta a^{-x}.$ Now, on using (\ref{eq1}), we get the required result. $\hfill\square$\\

\hspace*{.2in} Below we provide characterization of Type $3$ extreme value distribution.
\begin{t1}\label{th2.7}
For any absolutely continuous random variable $X$,
\begin{equation}\label{eq2.8}E[X\phi(X)]\geqslant\frac {2} {\eta^2-(1+c^2)},
\end{equation}
where $\eta=b/\mu$ and $c=\sigma/\mu.$
The equality holds if and only if $X$ follows the distribution given in (\ref{eq2}).
\end{t1}
Proof: By Cauchy-Schwartz inequality, we have
\begin{equation} \label{eq2.9} \left[\int_a^bxF(x)dx\right]\left[\int_a^b\frac {xf^2(x)}{F(x)} dx\right]\geqslant\mu^2.
\end{equation}
Since
$$\int_a^bxF(x)dx=\frac{b^2-(\sigma^2+\mu^2)}{2}$$
and
$$\int_a^b\frac {xf^2(x)}{F(x)}dx=E[X\phi(X)],$$
by noting $\eta=b/\mu(>0)$ and $c=\sigma/\mu(>0)$, (\ref{eq2.9}) reduces to (\ref{eq2.8}).
The equality in (\ref{eq2.9}) holds if and only if there exists some constant $A_0(>0)$ such that, for all $x\in(-\infty,b],$
$$xF(x)=A_0\frac {xf^2(x)}{F(x)}.$$
This gives $\phi(x)=1/\sqrt{A_0}~(constant)$. Hence the result follows from (\ref{eq1}).$\hfill\square$\\

\hspace*{.2in}In the following theorem we generalize the above result.
\begin{t1}\label{th2.8}
Let $X$ be an absolutely continuous random variable with $E[X^k\phi(X)]<\infty$. Then, for some integer $k>0$,
$$ E\left[X^k \phi(X)\right]\geqslant\frac {(k+1)\mu{_k}^{2}} {b^{k+1}-\mu_{k+1}}.$$
The equality holds if and only if $X$ belongs to the class of distribution given in (\ref{eq2}).
\end{t1}
Proof: By Cauchy-Schwartz inequality, we have
\begin{equation} \label{eq2.10} \left[\int_a^bx^kF(x)dx\right]\left[\int_a^b\frac {x^kf^2(x)}{F(x)}dx\right]\geqslant \mu{_k}^2.
\end{equation}
Clearly,
$$\int_a^bx^kF(x)dx=\frac{b^{k+1}-\mu_{k+1}}{k+1}$$
and
$$\int_a^b\frac {x^kf^2(x)}{F(x)}dx=E[X^k\phi(X)].$$
Now, from (\ref{eq2.10}), we get the required result.
The equality in (\ref{eq2.10}) holds if and only if there exists some constant $B_0(>0)$ such that, for all $x\in(-\infty,b],$
$$x^kF(x)=B_0\frac {x^kf^2(x)}{F(x)}.$$
This gives $\phi(x)=1/\sqrt{B_0}=\gamma,$ say. Hence the result follows from (\ref{eq1}).$\hfill\square$\\

\hspace*{.2in} A characterization of reflected Weibull distribution is given in the following theorem. For some properties and applications of reflected Weibull distribution in fitting sample data one may refer to Cohen (1973).
\begin{t1}\label{th2.9}
For any absolutely continuous random variable $X$,
$$ E\left[\frac{\phi(X)} {X}\right]\geqslant \frac {2} {b^2-(1+c^2) \mu^2},$$
and the equality holds if and only if $X$ is distributed as a subclass of reflected Weibull distribution defined on $(-\infty,0]$ and specified by
$$ F(x)=e^{-\theta x^2}~,~x\in(-\infty,0],~\theta>0.$$
\end{t1}
Proof: It is known that
$$\int_a^bxF(x)dx=\frac{b^2-(\sigma^2+\mu^2)}{2}.$$
By Cauchy-Schwartz inequality, we have
\begin{equation} \label{eq2.11} \left[\int_a^bxF(x)dx\right]\left[\int_a^b\frac {f^2(x)}{xF(x)}dx\right]\geqslant1,
\end{equation}
which is equivalent to the fact that
$$E\left[\frac{\phi(X)}{X}\right]\geqslant\frac {2} {b^2-\mu^2(1+c^2)}.$$
The equality in (\ref{eq2.11}) holds if and only if there exists some constant $\bar B(>0)$ such that, for all $x\in(-\infty,0],$
$$xF(x)=\bar B^2\frac {f^2(x)}{xF(x)}.$$
This gives $\phi(x)=-2\theta x$~where~$\theta(=1/2\bar B)>0$~is a constant. Now on using (\ref{eq1}) and noting the fact~$E(X)=\mu$, we get the required result.$\hfill\square$\\

\hspace*{.2in}In the following theorem we generalize the above result.
\begin{t1}\label{th2.10}
For any absolutely continuous random variable $X$ with $\mu_{k+1}<\infty$ and $E\left[\frac{\phi(X)} {X^k}\right]<\infty$,
$$ E\left[\frac{\phi(X)} {X^k}\right]\geqslant\frac{k+1}{b^{k+1}-\mu_{k+1}},$$
and the equality holds if and only if $X$ follows the  reflected Weibull distribution
$$ F(x)=e^{-\theta x^{k+1}}~;~x\in(-\infty,0],~\theta>0,$$
where k is odd positive integer.
\end{t1}
Proof: By Cauchy-Schwartz inequality, we have
\begin{equation} \label{eq2.12} \left[\int_a^bx^kF(x)dx\right]\left[\int_a^b\frac {f^2(x)}{x^kF(x)}dx\right]\geqslant1,
\end{equation}
which on simplification gives the required result.
The equality in (\ref{eq2.12}) holds if and only if there exists some constant $B^*(>0)$ such that, for all $x\in(-\infty,0],$
$$x^kF(x)=(B^*)^2\frac {f^2(x)}{x^kF(x)}.$$
This gives $\phi(x)=-\xi x$~where~$\xi(=1/B^*)>0$~is a constant. Now on using (\ref{eq1}) and noting the fact that $\theta=\xi/(k+1)(>0)$, we get the required result.
\section{Inequalities involving expectation of EIT}
Here we consider inequalities involving expectations of EIT to characterize quite a few continuous distributions used in reliability theory.
\begin{t1}\label{th3.1}
For any absolutely continuous random variable $X$,
\begin{equation} \label{eq3.1}E\left[\frac{1}{m(X)}\right]\geqslant\frac{1}{E[m(X)]},
\end{equation}
and equality holds if and only if $X$ follows the distribution given in (\ref{eq2}) with a change in the parameter.
\end{t1}
Proof: By Cauchy-Schwartz inequality, we obtain (\ref{eq3.1}).
The equality in (\ref{eq3.1}) holds if and only if there exists a constant $C(>0)$ such that, for all $x\in(-\infty,b]$,
$$\frac{f(x)}{m(x)}=C m(x)f(x),$$
which is equivalent to the fact that $m(x)=1/\sqrt{C}.$ Now, on using (\ref{eq0}) and Leibnitz rule of differentiation under the sign of integral, we get $\phi(x)=\sqrt{C}=\gamma$, say. Hence the result follows from (\ref{eq1}). $\hfill\square$\\

\hspace*{.2in} Below we consider another characterization result.
\begin{t1}\label{th3.2}
Let $X$ be an absolutely continuous nonnegative random variable with $E[X m(X)]<\infty$ and $E\left[\frac{1}{X m(X)}\right]<\infty $, then
\begin{equation} \label{eq3.2}E\left[\frac{1} {X m(X)}\right]\geqslant\frac{1}{E\left[Xm(X)\right]},
\end{equation}
and the equality holds if and only if $X$ follows the finite range distribution having cumulative distribution function
$$F(x)=(x/b)e^{\theta(x^2-b^2)}~,~\theta>0,~x\in(0,b].$$
\end{t1}
Proof: (\ref{eq3.2}) can easily be obtained from Cauchy-Schwartz inequality.
The equality in (\ref{eq3.2}) holds if and only if there exists a constant $C_0(>0)$ such that
$$\frac{f(x)}{xm(x)}=C_0 xm(x)f(x),$$
which is equivalent to the fact that $m(x)=\theta/x.$ Now, on using (\ref{eq0}), (\ref{eq1}) and Leibnitz rule of differentiation under the sign of integral, we get the required result. $\hfill\square$\\

\hspace*{.2in} In the following theorem we generalize the above result.
\begin{t1}\label{th3.3}
Let $X$ be an absolutely continuous nonnegative random variable with $E[X^k m(X)]<\infty$ and $E\left[\frac{1}{X^k m(X)}\right]<\infty $. Then, for some integer $k>0$,
\begin{equation} \label{eq3.3}E\left[\frac{1} {X^k m(X)}\right]\geqslant\frac{1}{E\left[X^k m(X)\right]},
\end{equation}
and the equality holds if and only if $X$ follows the finite range distribution having cumulative distribution function
$$F(x)=(x/b)^ke^{\theta(x^{k+1}-b^{k+1})},~\theta>0,~x\in(0,b].$$
\end{t1}
Proof: By Cauchy-Schwartz inequality, we obtain (\ref{eq3.3}).
The equality in (\ref{eq3.3}) holds if and only if there exists a constant $C'(>0)$ such that, for all $x\in(0,b]$,
$$\frac{f(x)}{x^k m(x)}=C' x^k m(x)f(x),$$
which is equivalent to the fact that $m(x)=\zeta x^{-k}$, where $\zeta(>0)$ is a constant. Now, on using (\ref{eq0}), (\ref{eq1}) and Leibnitz rule of differentiation under the sign of integral and also by noting $\theta=1/(\zeta(k+1))>0 $, we get the required result. $\hfill\square$\\

\hspace*{.2in} In the following two theorems we characterize some distributions given in Kundu et al. (2010).
\begin{t1}\label{th3.4}
For any absolutely continuous random variable $X$,
\begin{equation} \label{eq3.4}E\left[\frac{m(X)}{X}\right]\geqslant\frac{1}{E\left[\frac {X}{m(X)}\right]}.
\end{equation}
The equality holds if and only if $X$ follows the three distributions as given in Theorem 2.1 of Kundu et al. (2010) for suitable values of $p$.
\end{t1}
Proof: By Cauchy-Schwartz inequality, we obtain (\ref{eq3.4}).
The equality in (\ref{eq3.4}) holds if and only if there exists a constant $\bar C(>0)$ such that
$$\frac{m(x)}{x} f(x)=\bar C^2 \frac {x} {m(x)}f(x),$$
which is equivalent to the fact that $m(x)=\bar Cx$
Therefore, $\phi(x) m(x)=1-m^{'}(x)=1-\bar C=~p~(say).$
Thus, the result follows from Theorem 2.1 of Kundu et al. (2010). $\hfill\square$
\begin{t1}\label{th3.5}
For any absolutely continuous random variable $X$,
\begin{equation} \label{eq3.5}E\left[\frac{m(X)}{\alpha +\beta X}\right]\geqslant\frac{1}{E\left[\frac {\alpha +\beta X}{m(X)}\right]}.
\end{equation}
The equality holds if and only if $X$ follows the three distributions as given in Theorem 2.1 of Kundu et al. (2010) with  $p$ replaced by $(1-\xi\beta)$ for $\beta=0,~\beta>0$ and $\beta<0$ respectively where $\xi(>0)$ is a constant.
\end{t1}
Proof: (\ref{eq3.5}) is a simple consequence of Cauchy-Schwartz inequality.
The equality in (\ref{eq3.5}) holds if and only if there exists a constant $\xi(>0)$ such that
$$\frac{m(x)}{\alpha+\beta x} f(x)=\xi^2 \frac {\alpha+\beta x} {m(x)}f(x),$$
which is equivalent to the fact that $m(x)=\xi\beta x+\xi\alpha.$
Therefore, $\phi(x) m(x)=1-\xi\beta=~p~(say).$
Thus, the result follows from Theorem 2.1 of Kundu et al. (2010). $\hfill\square$\\

\hspace*{.2in} We conclude this section by giving two characterization results in analogy with Theorem \ref{th2.5} and Theorem \ref{th2.6}, respectively.
\begin{t1}\label{th3.6}
Let $X$ be an absolutely continuous random variable with $E[e^X m(X)]<\infty$ and $E\left[\frac{1}{e^X m(X)}\right]<\infty$, then
\begin{equation} \label{eq3.6}E\left[\frac{1} {e^X m(X)}\right]\geqslant\frac{1}{E\left[e^X m(X)\right]}
\end{equation}
and the equality holds if and only if $X$ follows a distribution having cumulative distribution function
$$F(x)=e^{(x-b)} e^{\theta(e^x-e^b)}~,~\theta>0,~x\in(-\infty,b].$$
\end{t1}
Proof: By Cauchy-Schwartz inequality, we obtain (\ref{eq3.6}).
The equality in (\ref{eq3.6}) holds if and only if there exists a constant $\theta(>0)$ such that, for all $x\in(-\infty,b]$,
$$\frac{f(x)}{e^xm(x)}=\theta^2{e^xm(x)f(x)},$$
which is equivalent to the fact that $m(x)=e^{-x}/\theta$. Now, on using (\ref{eq0}), (\ref{eq1}) and Leibnitz rule of differentiation under the sign of integral we get the required result. $\hfill\square$
\begin{t1}\label{th3.7}
Let $X$ be an absolutely continuous random variable with $E[a^X m(X)]<\infty$ and $E\left[\frac{1}{a^X m(X)}\right]<\infty$, where $a>1$. Then
\begin{equation} \label{eq3.7}E\left[\frac{1} {a^X m(X)}\right]\geqslant\frac{1}{E\left[a^X m(X)\right]}
\end{equation}
and the equality holds if and only if $X$ follows a distribution having cumulative distribution function
$$F(x)=e^{\gamma(x-b)} e^{\delta(a^x-a^b)}~,~\gamma>0,~\delta>0,~x\in(-\infty,b].$$
\end{t1}
Proof: The inequality in (\ref{eq3.7}) is a simple consequence of Cauchy-Schwartz inequality and
the equality holds if and only if there exists a constant $C^*(>0)$ such that, for all $x\in(-\infty,b]$,
$$\frac{f(x)}{e^xm(x)}={C^*}^2{e^xm(x)f(x)},$$
which is equivalent to the fact that $m(x)=a^{-x}/C^*$. Now, by noting that $\gamma=\log a(>0)$, $\delta=C^*/(\log a)(>0)$ and on using (\ref{eq0}), (\ref{eq1}) and Leibnitz rule of differentiation under the sign of integral we get the required result.
\section{Inequalities involving expectation of RHR with EIT and reversed AI function}
In this section we characterize Type 3 extreme value distribution through the expected values of RHR with EIT  and reversed AI function. Motivated by Theorem 2.5 of Nanda (2010) we have the following result.
\begin{t1}\label{th4.1}
For any absolutely continuous random variable $X$,
\begin{equation} \label{eq4.1}E\left[\frac{m(X)}{\phi (X)}\right]\geqslant\frac{1}{E\left[\frac {\phi (X)}{m(X)}\right]},
\end{equation}
and the equality holds if and only if $X$ follows the distribution given in (\ref{eq2}) for a suitable choice of $\gamma$.
\end{t1}
Proof: By Cauchy-Schwartz inequality, we obtain (\ref{eq4.1}). The equality in (\ref{eq4.1}) holds if and only if there exists
a constant $\xi(>0)$ such that for all $x\in(-\infty,b]$,
$$\frac{m(x)}{\phi (x)} f(x)=\xi^2 \frac {\phi (x)} {m(x)}f(x),$$
which is equivalent to the fact that
$$f(x)=\psi\int _{-\infty}^x F(u) du,~{\rm where}~\psi=1/\xi(>0).$$
Writing $\int _{-\infty} ^x F(u) du=y$, the above equation can be written as
$$ \frac {d^2y}{d^2x}=\psi y.$$
The solution of this second order differential equation is
$$\int _{-\infty}^x F(u)du=\gamma e^{\sqrt{\psi} x}+\delta e^{-\sqrt{\psi} x}.$$
Since $F(-\infty)=0$ and $F(b)=1$, on using $E(X)=\mu$, the above equation reduces to
$$ F(x)= e^{\frac {x-b}{b-\mu}},~x\in (-\infty,b].$$
Hence the result follows.$\hfill\square$\\

\hspace*{.2in} We generalize the above result in the following theorem.
\begin{t1}\label{th4.2}
Let $X$ be an absolutely continuous random variable. Then, for some integer $k\geqslant0$,
$$ E\left[X^k\frac{m(X)}{\phi (X)}\right]E\left[ X^k \frac {\phi (X)}{m(X)}\right]\geqslant{\mu_k}^2,$$
and the equality holds if and only if $X$ follows the distribution given in (\ref{eq2}).
\end{t1}
Proof: By Cauchy-Schwartz inequality, we have
\begin{equation} \label{eq4.2}E\left[X^k\frac{m(X)}{\phi (X)}\right]E\left[X^k\frac{\phi (X)}{m(X)}\right]\geqslant\mu_{k}^2.
\end{equation}
Thus the equality in (\ref{eq4.2}) holds if and only if there exists a constant $D(>0)$ such that for all $x\in(-\infty,b]$,
$$ x^k\frac{m(x)}{\phi (x)} f(x)=D^2x^k \frac {\phi (x)} {m(x)}f(x),$$
which is equivalent to the fact that
$$f(x)=\psi\int _{-\infty}^x F(u) du,$$
where $\psi=1/D(>0)$. Writing $\int _{-\infty} ^x F(u) du=y$, the above equation can be written as
$$ \frac {d^2y}{d^2x}=\psi y.$$
The solution of this second order differential equation is
$$\int _{-\infty}^x F(u)du=c_1 e^{\sqrt{\psi} x}+c_2 e^{-\sqrt{\psi} x}.$$
On using $F(-\infty)=0$ and $F(b)=1$, we get the required result.$\hfill\square$\\

\hspace*{.2in} In the sequel again we characterize Type 3 extreme value distribution through $E\left[\overline L(X)/\phi (X)\right]$.
\begin{t1}\label{th4.3}
For any absolutely continuous random variable $X$,
\begin{equation} \label{eq4.3}E\left[\frac{\overline L(X)}{\phi (X)}\right]\geqslant\frac{1}{E\left[\frac {\phi (X)}{\overline L(X)}\right]},
\end{equation}
and the equality holds if and only if $X$ follows the distribution given in (\ref{eq2}) for a different value of the parameter.
\end{t1}
Proof: By Cauchy-Schwartz inequality, we obtain (\ref{eq4.3}).
The equality in (\ref{eq4.3}) holds if and only if there exists a constant $\xi(>0)$ such that for all $x\in(-\infty,b]$,
$$\frac{\overline L(x)}{\phi (x)} f(x)=\xi^2 \frac {\phi (x)} {\overline L(x)}f(x).$$
Which is equivalent to the fact that
$$ \int_x^b\phi(u) du=\gamma (b-x),~{\rm where}~\gamma=1/\xi(>0)~{\rm is~a~constant}.$$
Differentiating above equation by Leibnitz rule of differentiation under the sign of integral, we get
$$ \phi(x)=\gamma,~{\rm a~constant}.$$
Hence the result follows from (\ref{eq1}).$\hfill\square$\\

\hspace*{.2in} A generalization of the above result is given in the following theorem.
\begin{t1}\label{th4.4}
Let $X$ be an absolutely continuous random variable with $E\left[X^k\frac{\overline L(X)}{\phi (X)}\right]<\infty$ and $E\left[ X^k \frac {\phi (X)}{\overline L(X)}\right]<\infty$. Then, for some integer $k\geqslant0$,
$$ E\left[X^k\frac{\overline L(X)}{\phi (X)}\right]E\left[ X^k \frac {\phi (X)}{\overline L(X)}\right]\geqslant~{\mu_k}^2,$$
and the equality holds if and only if $X$ follows the distribution given in (\ref{eq2}).
\end{t1}
Proof: By Cauchy-Schwartz inequality, we have
\begin{equation} \label{eq3.24}E\left[X^k\frac{\overline L(X)}{\phi (X)}\right]E\left[X^k\frac{\phi (X)}{\overline L(X)}\right]\geqslant~\mu_{k}^2 .
\end{equation}
Thus the equality in (\ref{eq3.24}) holds if and only if there exists a constant $D_0(>0)$ such that for all $x\in(-\infty,b]$,
$$ x^k\frac{\overline L(x)}{\phi (x)} f(x)=D_0^2x^k \frac {\phi (x)} {\overline L(x)}f(x).$$
Which is equivalent to the fact that
$$ \int_x^b\phi(u) du=\gamma (b-x),~{\rm where}~\gamma=1/D_0(>0)~{\rm is~a~constant}.$$
Hence the result follows from Theorem \ref{th4.3}.
\section*{Acknowledgements}
 The financial support (Ref. No. SR/FTP/MS-016/2012) from the Department of Science and Technology, Government of India is gratefully acknowledged by C. Kundu.

\end{document}